\documentclass[10pt]{article}

\usepackage{amscd,amstex,amssymb}

\def\eqref#1{(\ref{#1})}
\newcommand{\goth}{\mathfrak}
\newcommand{\arrow}{{\:\longrightarrow\:}}

\newcommand{\C}{{\Bbb C}}
\newcommand{\R}{{\Bbb R}}

\newcommand{\1}{\sqrt{-1}\:}
\newcommand{\inangles}[1]{{\langle #1\rangle}}
\newcommand{\galochka}{\check{\;}}
\newcommand{\restrict}[1]{{\left|_{{\phantom{|}\!\!}_{#1}}\right.}}

\renewcommand{\c}[1]{{\cal #1}}
\newcommand{\calo}{{\cal O}}

\renewcommand{\tilde}{\widetilde}
\renewcommand{\bar}{\overline}
\renewcommand{\phi}{\varphi}
\renewcommand{\epsilon}{\varepsilon}

\newcommand{\End}{\operatorname{End}}

\newcommand{\Def}{\operatorname{Def}}

\newcommand{\comment}[1]{{}}

\def\blacksquare{\hbox{\vrule width 4pt height 4pt depth 0pt}}
\def\endproof{\blacksquare}

\makeatletter

\@ifundefined{Bbb}
     {\newcommand{\Bbb}[1]{{\mathbb #1}}}%
{}

\newcommand{\ps@verbit}{%
  \renewcommand{\@oddhead}{%
          \scriptsize
          {Desingularization of hyperk\"ahler varieties I}
          \hfil\tiny {final version, November 7 1996}}
  \renewcommand{\@evenhead}{\@oddhead}
  \renewcommand{\@oddfoot}{\hfil\thepage\hfil}
  \renewcommand{\@evenfoot}{\@oddfoot}}
 
\newcounter{Mycounter}[section]
\newcounter{lemma}[section]
\renewcommand{\thelemma}{{Lemma \thesection.\arabic{lemma}}}
\newcommand{\lemma}{%
     \setcounter{lemma}{\value{Mycounter}}
     \refstepcounter{lemma}
     \stepcounter{Mycounter}
     {\bf \thelemma:\ }}

\newcounter{claim}[section]
\renewcommand{\theclaim}{{Claim \thesection.\arabic{claim}}}
\newcommand{\claim}{%
     \setcounter{claim}{\value{Mycounter}}
     \refstepcounter{claim}
     \stepcounter{Mycounter}
     {\bf \theclaim:\ }}

\newcounter{sublemma}[section]
\newcounter{corollary}[section]
\renewcommand{\thecorollary}{{Corollary \thesection.\arabic{corollary}}}
\newcommand{\corollary}{%
     \setcounter{corollary}{\value{Mycounter}}
     \refstepcounter{corollary}
     \stepcounter{Mycounter}
     {\bf \thecorollary:\ }}

\newcounter{theorem}[section]
\renewcommand{\thetheorem}{{Theorem \thesection.\arabic{theorem}}}
\newcommand{\theorem}{%
     \setcounter{theorem}{\value{Mycounter}}
     \refstepcounter{theorem}
     \stepcounter{Mycounter}
     {\bf \thetheorem:\ }}

\newcounter{conjecture}[section]
\renewcommand{\theconjecture}{{Conjecture \thesection.\arabic{conjecture}}}
\newcommand{\conjecture}{%
     \setcounter{conjecture}{\value{Mycounter}}
     \refstepcounter{conjecture}
     \stepcounter{Mycounter}
     {\bf \theconjecture:\ }}

\newcounter{proposition}[section]
\renewcommand{\theproposition}
       {{Proposition \thesection.\arabic{proposition}}}
\newcommand{\proposition}{%
     \setcounter{proposition}{\value{Mycounter}}
     \refstepcounter{proposition}
     \stepcounter{Mycounter}
     {\bf \theproposition:\ }}

\newcounter{definition}[section]
\renewcommand{\thedefinition}
       {{Definition \thesection.\arabic{definition}}}
\newcommand{\definition}{%
     \setcounter{definition}{\value{Mycounter}}
     \refstepcounter{definition}
     \stepcounter{Mycounter}
     {\bf \thedefinition:\ }}

\newcounter{example}[section]
\newcounter{remark}[section]
\renewcommand{\theremark}{{Remark \thesection.\arabic{remark}}}
\newcommand{\remark}{%
     \setcounter{remark}{\value{Mycounter}}
     \refstepcounter{remark}
     \stepcounter{Mycounter}
     {\bf \theremark:\ }}

\newcounter{problem}[section]
\newcounter{question}[section]
\@addtoreset{equation}{section}
\@addtoreset{footnote}{section}
\makeatother

\begin{document}

\begin{center}
{\Large\bf
	Desingularization of singular hyperk\"ahler varieties I.}\\[4mm]
Misha Verbitsky,\footnote{Supported by the NSF grant 9304580}\\[4mm]
{\tt verbit@@thelema.dnttm.rssi.ru, verbit@@math.ias.edu}
\end{center}

\hfill

{\small 
\hspace{0.2\linewidth}
\begin{minipage}[t]{0.7\linewidth}
Let $M$ be a singular hyperk\"ahler variety, obtained as a
moduli space of stable holomorphic bundles on a compact
hyperk\"ahler manifold (alg-geom/9307008). Consider $M$ as a complex
variety in one of the complex structures induced by the
hyperk\"ahler structure. We show that normalization of
$M$ is smooth, hyperk\"ahler and does not depend on
the choice of induced complex structure. 
\end{minipage}
}

\setcounter{section}{-1}
\section{Introduction}
\label{_Intro_Section_}


\setcounter{footnote}{1}

The structure of this paper is following.

\begin{itemize}

\item In the first section, we give a compendium 
of definitions and results from hyperk\"ahler geometry, all
known from literature.

\item Section \ref{_real_ana_Section_} deals 
with the real analytic
varieties underlying complex varieties. 
We define almost complex structures on a 
real analytic variety. This notion is used in
order to define hypercomplex varieties. We show
that a hyperk\"ahler manifold
is always hypercomplex.

\item In Section \ref{_singu_hype_Section_}, we 
give a definition of a singular hyperk\"ahler variety,
following \cite{_Verbitsky:Hyperholo_bundles_}  and 
\cite{_Verbitsky:Deforma_}. We cite basic properties
and list the examples of such manifolds.

\item In Section \ref{_SLHS_Section_}, we define locally homogeneous
singularities. A  space with locally homogeneous singularities 
(SLHS) is an analytic space $X$ such that for all $x\in X$, the
$x$-completion of a local ring $\calo_xX$ is isomorphic
to an $x$-completion of associated graded ring 
$(\calo_xX)_{gr}$. We show that a complex variety
is SLHS if and only if the underlying real analytic
variety is SLHS. This allows us to define invariantly 
the notion of a hyperk\"ahler SLHS. The natural examples
of hyperk\"ahler SLHS include the moduli spaces of stable
holomorphic bundles, considered in 
\cite{_Verbitsky:Hyperholo_bundles_}.
\footnote{In \cite{_Verbitsky:Hyperholo_bundles_},
we proved that the moduli of stable bundles over a 
compact hyperk\"ahler manifold is a hyperk\"ahler variety, if we
assume certain numerical restrictions on the bundle's 
Chern classes. The stable bundles satisfying these restrictions are
called {\bf hyperholomorphic}.}
We conjecture that every hyperk\"ahler variety is 
a space with locally homogeneous singularities.

\item In Section \ref{_tange_cone_Section_}, we study the 
tangent cone of a singular hyperk\"ahler manifold $M$
in the point $x\in M$. We show that its reduction,
which is a closed subvariety of $T_x M$, is a union 
of linear subspaces $L_i\subset T_x M$. These subspaces are
invariant under the natural quaternion action in
$T_x M$. This implies that a normalization of $(M,I)$
is smooth. Here, as usually,
$(M, I)$ denotes $M$ considered as a complex
variety, with $I$ a complex structure induced by the singular
hyperk\"ahler structure on $M$.

\item In Section \ref{_desingu_Section_}, we formulate and prove
the desingularization theorem for hyperk\"ahler varieties with locally
homogeneous singularities.
For each such variety $M$ we construct a finite
surjective morphism $\tilde M \stackrel n \arrow M$ 
of hyperk\"ahler varieties, such that $\tilde M$ is smooth
and $n$ is an isomorphism outside of singularities of $M$.
The $\tilde M$ is obtained as a normalization of $M$;
thus, our construction is canonical and functorial.


\end{itemize}


\section{Hyperk\"ahler manifolds}


\subsection{Definitions}

This subsection contains a compression of 
the basic definitions from hyperk\"ahler geometry, found, for instance, in
\cite{_Besse:Einst_Manifo_} or in \cite{_Beauville_}.
 
\hfill
 
\definition \label{_hyperkahler_manifold_Definition_} 
(\cite{_Besse:Einst_Manifo_}) A {\bf hyperk\"ahler manifold} is a
Riemannian manifold $M$ endowed with three complex structures $I$, $J$
and $K$, such that the following holds.
 
\begin{description}
\item[(i)]  the metric on $M$ is K\"ahler with respect to these complex 
structures and
 
\item[(ii)] $I$, $J$ and $K$, considered as  endomorphisms
of a real tangent bundle, satisfy the relation 
$I\circ J=-J\circ I = K$.
\end{description}
 
\hfill 
 
The notion of a hyperk\"ahler manifold was 
introduced by E. Calabi (\cite{_Calabi_}).

\hfill
 
Clearly, hyperk\"ahler manifold has the natural action of
quaternion algebra ${\Bbb H}$ in its real tangent bundle $TM$. 
Therefore its complex dimension is even.
For each quaternion $L\in \Bbb H$, $L^2=-1$,
the corresponding automorphism of $TM$ is an almost complex
structure. It is easy to check that this almost 
complex structure is integrable (\cite{_Besse:Einst_Manifo_}).
 
\hfill
 
\definition \label{_indu_comple_str_Definition_} 
Let $M$ be a hyperk\"ahler manifold, $L$ a quaternion satisfying
$L^2=-1$. The corresponding complex structure on $M$ is called
{\bf an induced complex structure}. The $M$ considered as a complex
manifold is denoted by $(M, L)$.
 
\hfill
 
Let $M$ be a hyperk\"ahler manifold. We identify the group $SU(2)$
with the group of unitary quaternions. This gives a canonical 
action of $SU(2)$ on the tangent bundle, and all its tensor
powers. In particular, we obtain a natural action of $SU(2)$
on the bundle of differential forms. 

\hfill

\lemma \label{_SU(2)_commu_Laplace_Lemma_}
The action of $SU(2)$ on differential forms commutes
with the Laplacian.
 
{\bf Proof:} This is Proposition 1.1
of \cite{_Verbitsky:Hyperholo_bundles_}. \endproof
 
Thus, for compact $M$, we may speak of the natural action of
$SU(2)$ in cohomology.

 
\subsection{Trianalytic subvarieties in compact hyperk\"ahler
manifolds.}
 
 
In this subsection, we give a definition and a few basic properties
of trianalytic subvarieties of hyperk\"ahler manifolds. 
We follow \cite{_Verbitsky:Symplectic_II_}.
 
\hfill
 
Let $M$ be a compact hyperk\"ahler manifold, $\dim_\R M =2m$.
 
\hfill
 
\definition\label{_trianalytic_Definition_} 
Let $N\subset M$ be a closed subset of $M$. Then $N$ is
called {\bf trianalytic} if $N$ is a complex analytic subset 
of $(M,L)$ for any induced complex structure $L$.
 
\hfill
 
Let $I$ be an induced complex structure on $M$,
and $N\subset(M,I)$ be
a closed analytic subvariety of $(M,I)$, $dim_\C N= n$.
Denote by $[N]\in H_{2n}(M)$ the homology class 
represented by $N$. Let $\inangles N\in H^{2m-2n}(M)$ denote 
the Poincare dual cohomology class. Recall that
the hyperk\"ahler structure induces the action of 
the group $SU(2)$ on the space $H^{2m-2n}(M)$.
 
\hfill
 
\theorem\label{_G_M_invariant_implies_trianalytic_Theorem_} 
Assume that $\inangles N\in  H^{2m-2n}(M)$ is invariant with respect
to the action of $SU(2)$ on $H^{2m-2n}(M)$. Then $N$ is trianalytic.
 
{\bf Proof:} This is Theorem 4.1 of 
\cite{_Verbitsky:Symplectic_II_}.
\endproof
 
\remark \label{_triana_dim_div_4_Remark_}
Trianalytic subvarieties have an action of quaternion algebra in
the tangent bundle. In particular,
the real dimension of such subvarieties is divisible by 4.


\subsection{Totally geodesic submanifolds.}


\nopagebreak
\hspace{5mm}
\proposition \label{_comple_geodesi_basi_Proposition_}
Let $X \stackrel \phi\hookrightarrow M$ be an embedding of Riemannian 
manifolds (not necessarily compact) compatible with the Riemannian
structure.
 Then the following conditions are equivalent.
 
\begin{description}
\item[(i)] Every geodesic line in $X$ is geodesic in $M$.
 
\item[(ii)] Consider the Levi-Civita connection $\nabla$ on $TM$,
and restriction of $\nabla$ to $TM \restrict{X}$. Consider the
orthogonal decomposition 
\begin{equation} \label{TM_decompo_Equation_} 
   TM\restrict{X} = TX \oplus TX^\bot. 
\end{equation}
Then, this decomposition is preserved by the connection $\nabla$.
\end{description}
 
{\bf Proof:} Well known; see, for instance, 
\cite{_Besse:Einst_Manifo_}.
\blacksquare
 
\hfill

\proposition \label{_triana_comple_geo_Proposition_} 
Let $X\subset M$ be a trianalytic submanifold of a hyperk\"ahler
manifold $M$, where $M$ is not necessarily compact. Then
$X$ is totally geodesic.

{\bf Proof:} This is \cite{_Verbitsky:Deforma_}, Corollary 5.4.
\endproof


\section{Real analytic varieties}
\label{_real_ana_Section_}


Let $X$ be a complex analytic variety. The ``real analytic
variety underlying $X$'' (denoted by $X_\R$)
is the following object. By definition, $X_\R$ is a ringed space
with the same topology as $X$, but with a different structure
sheaf, denoted by $\calo_{X_\R}$. Let $C(X, \R)$ be a sheaf
of continous $\R$-valued functions on $X$. Then $\calo_{X_\R}$
is a subsheaf of $C(X, \R)$, defined as follows. Let 
$A\subset C(X, \R)$ be an arbitrary subsheaf of $C(X, \R)$.
By $Ser(A)\subset C(X, \R)$, we denote the sheaf of all functions
which can be locally represented by the absolutely convergent
series $\sum P_i(a_1,..., a_n)$, where $a_1,..., a_n$ are
sections of $A$ and $P_i$ are polynomials with coefficients in $\R$.
By definition, $\calo_{X_\R}= Ser(Re \calo_X)$, where $Re \calo_X$
is a sheaf of real parts of holomorphic functions.
 
Another interesting sheaf associated with $X_\R$ is a sheaf
$\calo_{X_\R}\otimes\C\subset C(X,\C)$ of complex-valued
real analytic functions.

Consider the sheaf $\calo_X$ of holomorphic functions on $X$
as a subsheaf of the sheaf $C(X,\C)$ of continous $\C$-valued
functions on $X$. The sheaf $C(X,\C)$ has a natural authomorphism
$f\arrow \bar f$, where $\bar f$ is complex conjugation. 
By definition, the section $f$ of $C(X,\C)$  is called
{\bf antiholomorphic} if $\bar f$ is holomorphic.
Let $\calo_X$ be the sheaf of holomorphic functions,
and $\bar \calo_X$ be the sheaf of antiholomorphic
functions on $X$. Let $\calo_X \otimes_\C \bar\calo_X 
\stackrel i\arrow C(X, \C)$ be the natural multiplication 
map. Clearly, the image of $i$ belongs to the subsheaf
$\calo_{X_\R}\otimes \C\subset C(X, \C)$.
 
\hfill

\claim \label{_comple_real_ana_produ_Claim_}
The sheaf homomorphism $i:\; \calo_X \otimes_\C \bar\calo_X \arrow 
\calo_{X_\R}\otimes \C\subset C(X, \C)$
is injective. For each point $x\in X$, $i$ induces an isomorphism
on $x$-completions of $\calo_X \otimes_\C \bar\calo_X$
and $\calo_{X_\R}\otimes \C$.
 
{\bf Proof:} Well known (see, for instance,
\cite{_real_anal_spa:GMT_}). $\;\;\blacksquare$
 
\hfill

 Let $\Omega^1(\calo_{X_\R})$,
$\Omega^1(\calo_X \otimes_\C \bar\calo_X)$, $\Omega^1(\calo_{X_\R}\otimes \C)$
be the sheaves of K\"ahler differentials associated with the corresponding
ring sheaves. There are natural sheaf maps
 
\begin{equation} \label{_complexifi_of_Omega_Equation_}
   \Omega^1(\calo_{X_\R})\otimes \C \arrow 
   \Omega^1(\calo_{X_\R}\otimes \C)
\end{equation}
and 
\begin{equation} \label{_Omega_X_R_and_Omega_X_Equation_}
   \Omega^1(\calo_{X_\R}\otimes \C)\arrow 
   \Omega^1(\calo_X\otimes_\C \bar \calo_X),
\end{equation}
correspoding to the monomorphisms
 
\[ \calo_{X_\R}\hookrightarrow \calo_{X_\R}\otimes \C, \;\;
   \calo_X\otimes_\C \bar \calo_X\hookrightarrow\calo_{X_\R}\otimes \C
\]
 
\hfill

\claim \label{_differe_real_ana_and_co_ana_Claim_} 
The map \eqref{_complexifi_of_Omega_Equation_} is an isomorphism.
Tensoring both sides of \eqref{_Omega_X_R_and_Omega_X_Equation_} 
by $\calo_{X_\R}\otimes \C$ produces an isomorphism
\[ \Omega^1(\calo_X\otimes_\C \bar \calo_X) 
   \bigotimes_{\calo_X\otimes_\C \bar \calo_X}\bigg(\calo_{X_\R}\otimes \C\bigg)
   =\Omega^1(\calo_{X_\R}\otimes \C).
\]
 
{\bf Proof:} Clear. $\;\;\blacksquare$
 
\hfill

According to the general results about differentials
(see, for example, \cite{_Hartshorne:Alg_Geom_}, Chapter II,
Ex. 8.3), the sheaf
$\Omega^1(\calo_X\otimes_\C \bar \calo_X)$ admits a canonical
decomposition:
 
\[ \Omega^1(\calo_X\otimes_\C \bar \calo_X) =
   \Omega^1(\calo_X)\otimes_\C \bar \calo_X
   \oplus\calo_X\otimes_\C\Omega^1(\bar \calo_X).
\]
Let $\tilde I$ be an endomorphism of 
$\Omega^1(\calo_X\otimes_\C \bar \calo_X)$
which acts as a multiplication by $\1$ on
 
\[ \Omega^1(\calo_X)\otimes_\C \bar \calo_X
   \subset \Omega^1(\calo_X\otimes_\C \bar \calo_X)
\]
and as a multiplication by $-\1$ on
 
\[ \calo_X\otimes_\C\Omega^1(\bar \calo_X)
   \subset \Omega^1(\calo_X\otimes_\C \bar \calo_X).
\]
Let $\underline I$ be the corresponding 
$\calo_{X_\R}\otimes \C$-linear endomorphism of 
\[ \Omega^1(\calo_{X_\R})\otimes \C =
   \Omega^1(\calo_X\otimes_\C \bar \calo_X) 
   \otimes_{\calo_X\otimes_\C \bar \calo_X}
   \bigg(\calo_{X_\R}\otimes \C\bigg).
\]
As easy check ensures that $\underline I$
is {\it real}, that is, comes from the 
$\calo_{X_\R}$-linear endomorphism of $\Omega^1(\calo_{X_\R})$.
Denote this $\calo_{X_\R}$-linear endomorphism by
\[ 
   I:\; \Omega^1(\calo_{X_\R})\arrow \Omega^1(\calo_{X_\R}), 
\]
$I^2=-1$. The endomorphism $I$ is called {\bf a complex structure
operator}.  In the case when $X$ is smooth, $I$ coinsides with
the usual complex structure operator on the cotangent space.

\hfill
 
\definition\label{_commu_w_comple_str_Definition_} 
Let $X$, $Y$ be complex analytic varieties, and 
\[ f:\; X_\R\arrow Y_\R\] be a morphism of underlying real
analytic varieties. Let 
$f^* \Omega^1_{Y_\R} \stackrel P\arrow \Omega^1_{X_\R}$ be the
natural map of sheaves of differentials associated with $f$.
Let 
 
\[ I_X:\; \Omega^1_{X_\R}\arrow \Omega^1_{X_\R}, \;\;\;
   I_Y:\; \Omega^1_{Y_\R}\arrow \Omega^1_{Y_\R} 
\]
be the complex structure operators, and
\[ f^* I_Y:\; f^*\Omega^1_{Y_\R}\arrow f^*\Omega^1_{Y_\R} \]
be $\calo_{X_\R}$-linear automorphism of 
$f^*\Omega^1_{Y_\R}$ defined as a pullback of $I_Y$.
We say that $f$ {\bf commutes with the complex structure}
if 
 
\begin{equation}\label{_commu_w_comle_Equation_}
   P\circ f^* I_Y = I_X \circ P.
\end{equation}
 
\hfill

\theorem \label{_commu_w_comple_str_Theorem_} 
Let $X$, $Y$ be complex analytic varieties, and 
\[ f_\R:\; X_\R\arrow Y_\R\] be a morphism of underlying real
analytic varieties, which commutes with the complex structure.
Then there exist a morphism $f:\; X\arrow Y$ of 
complex analytic varieties, such that $f_\R$ 
is its underlying morphism.

\hfill

{\bf Proof:} By Corollary 9.4, \cite{_Verbitsky:Deforma_}, the map
$f$, defined on the sets of points of $X$ and $Y$,
is meromorphic; to prove \ref{_commu_w_comple_str_Theorem_},
we need to show it is holomorphic. Let $\Gamma \subset X \times Y$
be the graph of $f$. Since $f$ is meromorphic, $\Gamma$ is 
a complex subvariety of $X\times Y$.
It will suffice to show that the natural projections
$\pi_1:\; \Gamma \arrow X$, $\pi_2:\; \Gamma \arrow Y$ are
isomorphisms. By \cite{_Verbitsky:Deforma_}, Lemma 9.12, 
the morphisms $\pi_i$ are flat. Since $f_\R$ induces isomorphism
of Zariski tangent spaces, same is true of $\pi_i$. Thus,
$\pi_i$ are unramified. Therefore, the maps $\pi_i$ are
etale. Since they are one-to-one on
points, $\pi_i$ etale implies $\pi_i$ is an isomorphism.
 $\;\;\blacksquare$
 
\hfill

\definition 
Let $M$ be a real analytic variety, and
\[ I:\; \Omega^1(\calo_M)\arrow\Omega^1(\calo_M) \]
be an endomorphism satisfying $I^2=-1$. Then
$I$ is called {\bf an almost complex structure
on $M$}. If there exist a complex analytic structure $\goth C$
on $M$ such that $I$ appears as the complex structure operator
associated with $\goth C$, we say that $I$ is {\bf integrable}. 
\ref{_commu_w_comple_str_Theorem_} implies
that this complex structure is unique if it
exists. 
 
\hfill

\definition \label{_hypercomplex_Definition_}
(Hypercomplex variety)
Let $M$ be a real analytic variety equipped with almost
complex structures $I$, $J$ and $K$, such that
$I\circ J = -J \circ I = K$. Assume that for all
$a, b, c\in \R$, such that $a^2 + b^2 + c^2=1$,
the almost complex structure $a I + b J + c K$ is integrable.
Then $M$ is called {\bf a hypercomplex variety}.

\hfill

\noindent
\claim \label{_hyperka_hyperco_Claim_} 
Let $M$ be a hyperk\"ahler manifold. Then $M$ is hypercomplex.
{\bf Proof:} Let $I$, $J$ be induced complex structures.
We need to identify $(M, I)_\R$ and $(M,J)_\R$ in a natural way.
These varieties are canonically identified as $C^\infty$-manifolds;
we need only to show that this identification is real analytic.
This is \cite{_Verbitsky:Deforma_}, Proposition 6.5. \endproof

\hfill

The following proposition will be used further on in this paper.

\hfill

\proposition\label{_tange_cone_underly_Proposition_} 
Let $M$ be a complex variety, $x\in X$ a point, and $Z_xM\subset T_xM$
be the reduction of the Zariski tangent cone to $M$ in $x$, considered
as a closed subvariety of the Zariski tangent space $T_xM$.
Let $Z_x M_\R \subset T_x M_\R$ be the Zariski tangent cone for the
underlying real analytic variety $M_\R$. Then 
$(Z_x M)_\R \subset (T_x M)_\R = T_x M_\R$ coinsides with $Z_x M_\R$.

{\bf Proof:} For each $v\in T_x M$, the point $v$ belongs to
$Z_x M$ if and only if there exist a real analytic path
$\gamma:\; [0, 1] \arrow M$, $\gamma(0)=x$ satisfying 
$\frac{d\gamma}{dt}=v$. The same holds true for $Z_x M_\R$.
Thus, $v\in Z_x M$ if and only if $v\in Z_x M_\R$. \endproof


\section{Singular hyperk\"ahler varieties.}
\label{_singu_hype_Section_}


In this section, we follow \cite{_Verbitsky:Deforma_},
Section 10. For more examples, motivations and reference,
the reader is advised to check \cite{_Verbitsky:Deforma_}.

\hfill

\definition\label{_singu_hype_Definition_}
(\cite{_Verbitsky:Hyperholo_bundles_}, Definition 6.5)
Let $M$ be a hypercomplex variety (\ref{_hypercomplex_Definition_}).
The following data define a structure of {\bf hyperk\"ahler variety}
on $M$.

\begin{description}

\item[(i)] For every $x\in M$, we have an $\R$-linear 
symmetric positively defined
bilinear form $s_x:\; T_x M \times T_x M \arrow \R$
on the corresponding real Zariski tangent space.

\item[(ii)] For each triple of induced complex structures
$I$, $J$, $K$, such that $I\circ J = K$, we have a 
holomorphic differential 2-form $\Omega\in \Omega^2(M, I)$.

\item[(iii)] 
Fix a triple of induced complex structure 
$I$, $J$, $K$, such that $I\circ J = K$. Consider the
corresponding differential 2-form $\Omega$ of (ii).
Let $J:\; T_x M \arrow T_x M$ be an endomorphism of 
the real Zariski tangent spaces defined by $J$, and $Re\Omega\restrict x$
the real part of $\Omega$, considered as a bilinear form on $T_x M$.
Let $r_x$ be a bilinear form $r_x:\; T_x M \times T_x M \arrow \R$ 
defined by $r_x(a,b) = - Re\Omega\restrict x (a, J(b))$.
Then $r_x$ is equal to the form $s_x$ of (i). In particular,
$r_x$ is independent from the choice of $I$, $J$, $K$.

\end{description}

\noindent \remark \label{_singu_hype_Remark_}
\nopagebreak
\begin{description}
\item[(a)] It is clear how to define a morphism of hyperk\"ahler varieties.

\item[(b)]
For $M$ non-singular,  \ref{_singu_hype_Definition_} is
 equivalent to the usual
one (\ref{_hyperkahler_manifold_Definition_}). 
If $M$ is non-singular,
the form $s_x$ becomes the usual Riemann form, and 
$\Omega$ becomes the standard holomorphically symplectic form.

\item[(c)] It is easy to check the following.
Let $X$ be a hypercomplex subvariety of a hyperk\"ahler
variety $M$. Then, restricting the forms $s_x$ and $\Omega$
to $X$, we obtain a hyperk\"ahler structure on $X$. In particular,
trianalytic subvarieties of hyperk\"ahler manifolds are always
hyperk\"ahler, in the sense of \ref{_singu_hype_Definition_}.

\end{description}

\hfill

{\bf Caution:} Not everything which is seemingly hyperk\"ahler
satisfies the conditions of \ref{_singu_hype_Definition_}.
Take a quotient $M/G$ os a hyperk\"ahler manifold by an action 
of finite group $G$, acting in accordance with hyperk\"ahler
structure. Then $M/G$ is, generally speaking, {\it not} hyperk\"ahler
(see \cite{_Verbitsky:Deforma_}, Section 10 for details).

\hfill

The following theorem, proven in 
\cite{_Verbitsky:Hyperholo_bundles_} (Theorem 6.3), 
gives a convenient way to construct
examples of hyperk\"ahler varieties.

\hfill

\theorem \label{_hyperho_defo_hyperka_Theorem_}
Let $M$ be a compact hyperk\"ahler manifold, $I$ an induced
complex structure and $B$ a stable holomorphic bundle over $(M, I)$.
Let $\Def(B)$ be a reduction\footnote{The deformation space might have
nilpotents in the structure sheaf. We take its reduction to avoid
this.} of the
deformation space of stable holomorphic structures on $B$.
Assume that $c_1(B)$, $c_2(B)$ are $SU(2)$-invariant, with respect
to the standard action of $SU(2)$ on $H^*(M)$. Then $\Def(B)$ has a
natural structure of a hyperk\"ahler variety. 

\nopagebreak
\endproof


\section{Spaces with locally homogeneous singularities}
\label{_SLHS_Section_}

 
\noindent
\definition
(local rings with LHS)
Let $A$ be a local ring. Denote by $\goth m$ its maximal ideal.
Let $A_{gr}$ be the corresponding associated graded ring.
Let $\hat A$, $\hat A_{gr}$ be the $\goth m$-adic completion
of $A$, $A_{gr}$. Let $(\hat A)_{gr}$, $(\hat A_{gr})_{gr}$ 
be the associated graded rings, which are naturally isomorphic to
$A_{gr}$. We say that $A$ {\bf has locally homogeneous singularities}
(LHS)
if there exists an isomorphism $\rho:\; \hat A \arrow \hat A_{gr}$
which induces the standard isomorphism 
$i:\; (\hat A)_{gr}\arrow (\hat A_{gr})_{gr}$ on associated
graded rings.

\hfill

\definition
(SLHS)
Let $X$ be a complex or real analytic space. Then 
$X$ is called {be a space with locally homogeneous singularities}
(SLHS) if for each $x\in M$, the local ring $\calo_x M$ 
has locally homogeneous singularities.

\hfill

By {\bf system of coordinates} on a complex space $X$, defined
in a neighbourhood $U$ of $x\in X$, we understand
a closed embedding $U\hookrightarrow B$ where $B$ is an open
subset of $\C^n$. Clearly, a system of coordinates can be considered as
a set of functions $f_1, ..., f_n$ on $U$. Then $U \subset B$
is defined by a system of equations on $f_1, ... f_n$.

\hfill

\remark \label{_SLHS_term_expla_Remark_}
Let $X$ be a complex space. Assume that for each $x\in X$, there 
exist a system of coordinates $f_1, ... , f_n$ in a neighbourhood 
$U$ of $x$, such that $U\subset B$ is defined by a system of homogeneous
polynomial equations. Then $X$ is a 
space with locally homogeneous singularities.
This explains the term.

\hfill

\claim \label{_redu_SLHS_Claim_}
Let $X$ be a complex analytic space with 
locally homogeneous singularities, and $X_r$ its reduction
(same space, with structure sheaf factorized by nilradical).
Then $X_r$ is also a space with 
locally homogeneous singularities.

{\bf Proof:} Clear. \endproof

\hfill

\lemma \label{_produ_rings_LHS_Lemma_} 
Let $A_1$, $A_2$ be local rings over $\C$, with
$A_i/ {\goth m}_i = \C$, where ${\goth m}_i$ is the
maximal ideal of $A_i$. Then $A_1\otimes_\C A_2$ is
LHS if and only if $A_1$ and $A_2$ are LHS.

\hfill

{\bf Proof (``if'' part):} Let 
$\rho_i:\; \hat A_i \arrow \widehat {(A_i)_{gr}}$ be the maps 
given by LHS condition.  Consider the map
\begin{equation} \label{_produ_of_LHS_Equation_}
   \rho_1 \otimes \rho_2:\; \hat A_i \otimes_\C \hat A_2
   \arrow \widehat{(A_i)_{gr}} \otimes_\C \widehat{(A_2)_{gr}}.
\end{equation}
Denote the functor of adic completions of local rings by
$B \arrow \widehat B$. Clearly, 
$\widehat{\hat A_i \otimes_\C \hat A_2} =
\widehat{A_1\otimes_\C A_2}$, and 
$\widehat{(\hat A_i)_{gr} \otimes_\C (\hat A_2)_{gr}} =
\widehat{(A_1)_{gr}\otimes_\C (A_2)_{gr}}$.
Plugging these isomorphisms into the completion
of both sides of
\eqref{_produ_of_LHS_Equation_}, we obtain that
a completion of $\rho_1 \otimes \rho_2$ provides an
LHS map for $A_1\otimes_\C A_2$. 

\hfill

{\bf ``only if'' part:} Let 
\[ \rho:\; \widehat{A_1\otimes_\C A_2} \arrow 
\widehat{((A_1)\otimes_\C (A_2))_{gr}}\] be the LHS map for 
$A_1\otimes_\C A_2$. There are natural maps 
\[u:\; \hat A_1 \arrow \widehat{A_1\otimes_\C A_2} \]
and \[ v:\; \widehat{((A_1)\otimes_\C (A_2))_{gr}} \arrow
(\hat A_1)_{gr}.\] The $u$ comes from the natural embedding
$a \arrow a\otimes 1\in A_1\otimes_\C A_2$ and $v$ from
the natural surjection 
$a\otimes b \arrow a \otimes \pi(b) \in A_1\otimes_\C \C$,
where $\pi:\; A_2 \arrow \C$ is the standard quotient map.
It is clear that $u\circ v$ induces identity on the associated graded
ring of $A_1$. \ref{_produ_rings_LHS_Lemma_} is proven.
\endproof

\hfill

\proposition\label{_comple_LHS<=>real_LHS_Proposition_}
Let $M$ be a complex variety, $M_\R$ the underlying real analytic
variety. Then $M_\R$ is a space with locally
homogeneous singularities (SLHS) if and only if $M$
is SLHS.

{\bf Proof:} By \ref{_comple_real_ana_produ_Claim_},
$\widehat{(\calo_x M_\R)\otimes \C} = 
\widehat{\calo_x M \otimes \bar\calo_x M}$. Thus,
\ref{_comple_LHS<=>real_LHS_Proposition_} is implied immediately
by \ref{_produ_rings_LHS_Lemma_}. \endproof

\hfill

\corollary \label{_hype_SLHS_for_diff_indu_c_str_Corollary_}
Let $M$ be a hyperk\"ahler variety, $I_1$, $I_2$ induced complex
structures. Then $(M, I_1)$ is a space with locally
homogeneous singularities if and only is $(M, I_2)$ is 
SLHS.

{\bf Proof:} The real analytic variety underlying 
 $(M, I_1)$ coinsides with that underlying 
 $(M, I_2)$. Applying \ref{_comple_LHS<=>real_LHS_Proposition_},
we immediately 
obtain \ref{_hype_SLHS_for_diff_indu_c_str_Corollary_}.
\endproof

\hfill

\definition
Let $M$ be a hyperk\"ahler variety. Then $M$ is called 
a space with locally homogeneous singularities (SLHS) if 
the underlying real analytic variety is SLHS
or, equivalently, for some induced complex structure
$I$ the $(M, I)$ is SLHS.

\hfill

\theorem 
Let $M$ be a compact hyperk\"ahler manifold, $I$ an induced
complex structure and $B$ a stable holomorphic bundle over $(M, I)$.
Assume that $c_1(B)$, $c_2(B)$ are $SU(2)$-invariant, with respect
to the standard action of $SU(2)$ on $H^*(M)$. 
Let $\Def(B)$ be a reduction
of a deformation space of stable holomorphic 
structures on $B$, which is a hyperk\"ahler variety by 
\ref{_hyperho_defo_hyperka_Theorem_}. Then 
$\Def(B)$ is a space with locally homogeneous singularities (SLHS).

\hfill

{\bf Proof:} Let $x$ be a point of $\Def(B)$, corresponding
to a stable holomorphic bundle $B$. In
\cite{_Verbitsky:Hyperholo_bundles_}, Section 7,
the neighbourhood $U$
of $x$ in $\Def(B)$ was described explicitely as follows. We constructed
a locally closed holomorphic embedding 
$U\stackrel \phi\hookrightarrow H^1(\End(B))$.
We proved that $v\in H^1(\End(B))$ belongs to the image of 
$\phi$ if and only if $v^2 =0$. Here $v^2 \in H^2(\End(B))$
is the square of $v$, taken with respect to the product
\[ H^1(\End(B))\times H^1(\End(B)) \arrow  H^2(\End(B)) \]
associated with the algebraic structure on $\End(B)$.
Clearly, the relation $v^2 =0$ is homogeneous. This relation 
defines a locally closed SLHS subspace $Y$ of $H^1(\End(B))$,
such that $\phi(U)$ is its reduction. Applying
\ref{_redu_SLHS_Claim_}, we obtain that $\phi(U)$ is also
a space with locally homogeneous singularities.
\endproof

\hfill

\conjecture \label{_hype_SLHS_Conjecture_}
Let $M$ be a hyperk\"ahler variety. Then $M$ is 
a space with locally homogeneous singularities.

\hfill

There is a rather convoluted argument which might prove
\ref{_hype_SLHS_Conjecture_}. This argument will be a subject 
of forthcoming paper \cite{_Verbitsky:singuII_}.


\section{Tangent cone of a hyperk\"ahler variety}
\label{_tange_cone_Section_}


Let $M$ be a hyperk\"ahler variety, $I$ an induced complex
structure and $Z_x(M,I)$ be a reduction of 
a Zariski tangent cone to $(M,I)$ in $x\in M$. 
Consider $Z_x(M,I)$ as a closed subvariety in 
the Zariski tangent space $T_xM$. The space $T_xM$ has a natural
metric and quaternionic structure. This makes $T_xM$ into a 
hyperk\"ahler manifold, isomorphic to ${\Bbb H}^n$. 

\hfill

\noindent\theorem \label{_cone_hype_Theorem_}
Under these assumptions, the following assertions hold:
\begin{description}
\item[(i)] The subvariety $Z_x(M,I)\subset T_x M$ is independent
from the choice of induced complex structure $I$.
\item [(ii)] Moreover, $Z_x(M,I)$ is a trianalytic subvariety
of $T_x M$.
\end{description}

{\bf Proof:} \ref{_cone_hype_Theorem_} (i) is implied by
\ref{_tange_cone_underly_Proposition_}. 
By \ref{_cone_hype_Theorem_} (i),
the Zariski tangent cone $Z_x(M,I)$ is a hypercomplex subvariety of $TM$.
According to \ref{_singu_hype_Remark_} (c), this implies that
$Z_x(M)$ is hyperk\"ahler. \endproof

\hfill

Further on, we denote the Zariski tangent cone to a hyperk\"ahler variety
by $Z_xM$. The Zariski tangent cone is equipped with a natural
hyperk\"ahler structure. 

\hfill

The following theorem shows that the Zariski tangent cone
$Z_xM\subset T_x M$ is a union of planes $L_i\subset T_x M$.

\hfill

\theorem \label{_cone_flat_Theorem_}
Let $M$ be a hyperk\"ahler variety, $I$ an induced complex
structure and $x\in M$ a point. 
Consider the reduction of the
Zariski tangent cone (denoted by $Z_x M$) as a subvariety of the
quaternionic space $T_x M$. Let $Z_x(M, I)= \cup L_i$ 
be the irreducible decomposition of the complex variety $Z_x(M,I)$.
Then

\begin{description}
\item[(i)] The decomposition $Z_x(M, I)= \cup L_i$ 
is independent from the choice of induced complex structure $I$.
\item[(ii)] For every $i$, the variety 
$L_i$ is a linear subspace of $T_x M$,
invariant under quaternion action.
\end{description}

{\bf Proof:}  Let $L_i$ be an irreducible component of
$Z_x(M, I)$, $Z_x^{ns}(M,I)$ be the non-singular part
of $Z_x(M,I)$, and $L_i^{ns}:=Z_x^{ns}(M,I) \cap L_i$.
Then $L_i$ is a closure of $L_i^{ns}$ in $T_xM$.
Clearly from \ref{_cone_hype_Theorem_}, $L_i^{ns}(M)$ is a hyperk\"ahler
submanifold in $T_xM$. By 
\ref{_triana_comple_geo_Proposition_}, $L_i^{ns}$ is totally
geodesic. A totally geodesic submanifold of a flat manifold is
again flat. Therefore, $L_i^{ns}$ is an open subset of a linear
subspace $\tilde L_i\subset T_xM$. Since $L_i^{ns}$ is a hyperk\"ahler
submanifold, $\tilde L_i$ is invariant with respect to quaternions.
The closure $L_i$ of $L_i^{ns}$ is a complex analytic subvariety
of $T_x(M,I)$. Therefore, $\tilde L_i = L_i$. This proves 
\ref{_cone_flat_Theorem_} (ii). From the above argument, it is
clear that $Z_x^{ns}(M,I)= \coprod L_i^{ns}$ (disconnected sum).
Taking connected components of $Z_x^{ns}M$ for each induced 
complex structure, we obtain the same decomposition
$Z_x(M, I)= \cup L_i$, with $L_i$ being closured of connected components.
This proves \ref{_cone_flat_Theorem_} (ii). \endproof

\hfill

\corollary\label{_normali_smooth_Corollary_}
Let $M$ be a hyperk\"ahler variety, and $I$ an induced 
complex structure. Assume that $M$ is a 
space with locally homogeneous singularities.
Then the normalization of $(M,I)$ is smooth.

{\bf Proof:} The normalization of $Z_xM$ is smooth by 
\ref{_cone_flat_Theorem_}. The 
normalization is compatible with the adic
completions (\cite{_Matsumura:Commu_Alge_},
Chapter 9, Proposition 24.E). Therefore, the integral closure
of the completion of $\calo_{Z_xM}$ is a regular ring.
Now, from the definition of locally homogeneous 
intersections, it follows that the integral closure of
$\calo_xM\galochka$ is also a regular ring, where 
$\calo_xM\galochka$ is an adic completion of the
local ring of holomorphic functions on $(M, I)$ in a neighbourhood
of $x$. Applying \cite{_Matsumura:Commu_Alge_},
Chapter 9, Proposition 24.E again, we obtain that the
integral closure of $\calo_x M$ is regular. This proves
\ref{_normali_smooth_Corollary_}
\endproof


\section{Desingularization of hyperk\"ahler varieties}
\label{_desingu_Section_}


\noindent\theorem \label{_desingu_Theorem_}
Let $M$ be a hyperk\"ahrler or a hypercomplex variety. 
Assume that $M$ is a space with 
locally homogeneous singularities, and 
$I$ an induced complex structure.
Let \[ \widetilde{(M, I)}\stackrel n\arrow (M,I)\] 
be the normalization of
$(M,I)$. Then $\widetilde{(M, I)}$ is smooth and
has a natural hyperk\"ahler structure $\c H$, such that the associated
map $n:\; \widetilde{(M, I)} \arrow (M,I)$ agrees with $\c H$.
Moreover, the hyperk\"ahler manifold $\tilde M:= \widetilde{(M, I)}$
is independent from the choice of induced complex structure $I$.

\hfill

{\bf Proof:} The variety $\widetilde{(M, I)}$
is smooth by \ref{_normali_smooth_Corollary_}.
Let $x\in M$, and $U\subset M$ be a neighbourhood of $x$.
Let ${\goth R}_x(U)$ be the set of irreducible components of 
$U$ which contain $x$. There is a natural map
$\tau: {\goth R}_x(U) \arrow Irr(Spec \calo_xM\galochka)$,
where $Irr(Spec \calo_xM\galochka)$ is a set of irreducible 
components of $Spec \calo_xM\galochka$, where 
$\calo_xM\galochka$ is a completion of 
$\calo_xM$ in $x$. Since $\calo_x M$ is Henselian 
(\cite{_Raynaud_}, VII.4), there exist a neighbourhood $U$ of $x$
such that $\tau: {\goth R}_x(U) \arrow Irr(Spec \calo_xM\galochka)$
is a bijection. Fix such an $U$. Since $M$ is a space locally
with locally homogeneous singularities, the irreducible decomposition
of $U$ coinsides with the irreducible decomposition of
the tangent cone $Z_x M$. 

Let $\coprod U_i \stackrel u \arrow U$ be the morphism
mapping a disjoint union of irreducible components of $U$ 
to $U$. By \ref{_cone_flat_Theorem_}, the $x$-completion
of $\calo_{U_i}$ is regular. Shrinking $U_i$ if necessary,
we may assume that $U_i$ is smooth. Then, the morphism
$u$ coinsides with the normalization of $U$. 

For each variety $X$, we denote by $X^{ns}\subset X$ 
the set of non-singular
points of $X$. Clearly, $u(U_i) \cap U^{ns}$ is a connected
component of $U^{ns}$. Therefore, $u(U_i)$ is trianalytic
in $U$. By \ref{_singu_hype_Remark_} (c), $U_i$ has a natural
hyperk\"ahler structure, which agrees with the map $u$.
This gives a hyperk\"ahler structure on the normalization
$\tilde U := \coprod U_i$. Gluing these hyperk\"ahler structures,
we obtain a hyperk\"ahler structure $\c H$ on 
the smooth manifold$\widetilde{(M, I)}$.
Consider the normalization map $n:\; \widetilde{(M, I)} \arrow M$,
and let $\tilde M^{n}:= n^{-1}(M^{ns})$. Then, 
$n\restrict{\tilde M^{n}} \tilde M^{n}\arrow M^{ns}$ 
is a finite covering which is compatible with the hyperk\"ahler
structure. Thus, $\c H\restrict{\tilde M^{n}}$ can be obtained
as a pullback from $M$. Clearly, a hyperk\"ahler structure 
on a manifold is uniquely
defined by its restriction to an open dense subset. We
obtain that $\c H$ is independent from the choice of $I$.
\endproof

\hfill

\remark
The desingularization argument works well for hypercomplex varieties.
The word ``hyperk\"ahler'' in this article can be
in most cases replaced by ``hypercomplex'', because we never
use the metric structure.

\hfill

{\bf Acknowledgements:} It is a pleasure to acknowledge 
the help of P. Deligne, who pointed out an error in the original argument.
Deligne also suggested the term ``locally homogeneous singularities''.
I am grateful to A. Beilinson, D. Kaledin, D. Kazhdan, T. Pantev
and S.-T. Yau for enlightening discussions.


\begin{thebibliography}{6666}
 


\bibitem[Beau]{_Beauville_} 
Beauville, A. Varietes K\"ahleriennes dont la premiere classe de Chern est
nulle. // J. Diff. Geom. 18, p. 755-782 (1983).


\bibitem[Bes]{_Besse:Einst_Manifo_} 
Besse, A., {\em Einstein Manifolds}, Springer-Verlag, New York (1987)

\bibitem[C]{_Calabi_} Calabi,  E.,
{\em Metriques k\"ahleriennes et fibr\`es holomorphes}, 
Ann. Ecol. Norm. Sup. {\bf 12} (1979), 269-294.  

\bibitem[GMT]{_real_anal_spa:GMT_} Guaraldo, F., Macri, P., Tancredi, A.,
{\em Topics on real analytic spaces},
 Advanced lectures in mathematics, Braunschweig: F. Vieweg, 1986.

\bibitem[H]{_Hartshorne:Alg_Geom_} 
 Hartshorne, R., {\em Algebraic geometry}, //
 Graduate texts in mathematics, vol. 52,
 New York : Springer-Verlag, 1977.
 

\bibitem[M] {_Matsumura:Commu_Alge_}
 Hideyuki Matsumura, Commutative algebra. // 
 Mathematics lecture note series, vol. 56,
 Benjamin/Cummings Pub. Co., Reading, Mass.,  1980.
 


\bibitem[R]{_Raynaud_} Raynaud, M. {\em Anneaux Locaux Hens\'eliens}, 
Springer LNM 169, 1970.


\bibitem[V-bun]{_Verbitsky:Hyperholo_bundles_} 
Verbitsky M., 
{\em Hyperholomorphic bundles over a hyperk\"ahler manifold}, 
alg-geom electronic preprint 9307008 (1993), 43 pages, LaTeX,\\
also published in: 
Journ. of Alg. Geom., {\bf 5} no. 4 (1996) pp. 633-669.


\bibitem[V2]{_Verbitsky:Symplectic_II_} 
Verbitsky M., {\em Hyperk\"ahler embeddings and holomorphic 
symplectic geometry II,} alg-geom electronic preprint 9403006 (1994),
14 pages, LaTeX,
also published in: GAFA {\bf 5} no. 1 (1995), 92-104.

\bibitem[V3]{_Verbitsky:Deforma_}
Verbitsky M., {\em Deformations of trianalytic subvarieties of
hyperk\"ahler manifolds}, alg-geom electronic preprint 9610010
(1996), 51 pages, LaTeX2e.

\bibitem[V-ne]{_Verbitsky:singuII_}
Verbitsky M.,
{\em Desingularization of singular hyperk\"ahler varieties II,}
(tentative), forthcoming.


\end{thebibliography}
\end{document}